\documentclass[12pt,reqno]{article}
\usepackage{amsmath,amsthm,amssymb,amscd}

\theoremstyle{plain} %% This is the default
\newtheorem{theorem}{Theorem}

\theoremstyle{definition}
\newtheorem{definition}{Definition}

\theoremstyle{remark}

\numberwithin{equation}{section}

\newcommand{\secref}[1]{Section~\ref{#1}}

\newcommand{\dem}[1]{{\em{#1}}}  % Emphasis within definitions

\newcommand{\ot}{\otimes}
\renewcommand{\t}{\tilde}

\newcommand{\lie}{\mathfrak{g}}	% Lie algebra
\newcommand{\qlie}{\mathfrak{g}_h}	% abstract quantum Lie algebra
\newcommand{\lqlie}{\mathfrak{L}_h(\lie)} % quantum Lie algebra inside
\newcommand{\uqg}{U_h(\lie)}   % quantized enveloping algebra
\newcommand{\sln}{\mathfrak{sl}_{n}}
\newcommand{\slt}{\mathfrak{sl}_{2}}
\newcommand{\lqslt}{\mathfrak{L}_h(\slt)} % quantum sl_2
\newcommand{\qslt}{(\slt)_h}
\newcommand{\ch}{{\mathbb C}[[h]]}
\newcommand{\CC}{{\mathbb C}}
\newcommand{\NN}{{\mathbb N}}
\newcommand{\lb}{[,]}		% Lie bracket
\newcommand{\ad}[1]{(\text{ad}\,#1)\,}   % adjoint action
\newcommand{\xp}{X^+_h}
\newcommand{\xm}{X^-_h}
\newcommand{\xpm}{X^\pm_h}
\newcommand{\hh}{H_h}
\newcommand{\poly}[1]{{\cal{P}}(#1)}
\newcommand{\qconj}{\sim}
\renewcommand{\t}[1]{\tilde{#1}}
\newcommand{\mat}{\text{Mat}}

\begin{document}

\title{
\vspace{-15mm}
\begin{flushright}\small q-alg/9605026 \\[30pt]\end{flushright}
\Huge Introduction to\\
\Huge Quantum Lie Algebras\footnote{
Contribution to the Proceedings of the Banach Minisemester on Quantum
Groups, Warsaw, November 1995}}

\author{Gustav W. Delius\footnote{
On leave from Department of Physics, Bielefeld University, Germany}
\\[8pt] 
  \small Department of Mathematics, King's College London\\[-4pt]
  \small Strand, London WC2R 2LS, Great Britain\\[-4pt]
  \small e-mail: delius@mth.kcl.ac.uk\\[-5pt]
  \small http://www.mth.kcl.ac.uk/\~{}delius/}
\vspace{4mm} 

\date{}
\maketitle

\begin{abstract}
Quantum Lie algebras are generalizations of Lie algebras whose
structure constants are power series in $h$.
They are derived from the quantized enveloping algebras $\uqg$.
The quantum Lie bracket satisfies a generalization of antisymmetry.
Representations of quantum Lie algebras are defined in terms of
a generalized commutator.

In this paper the recent general results about quantum Lie algebras are introduced
with the help of the explicit example of $(sl_2)_h$.
\end{abstract}

\thispagestyle{empty}
\newpage

\section{Introduction}

The subject of this paper are the question marks in the following diagram
\[
%\text{classical}~~~~
%\mbox{
\begin{CD}
\lie @>?>> \qlie\\
@VVV @AA?A\\
U(\lie) @>\text{Drinfel'd}>\text{Jimbo}> U_h(\lie)
\end{CD}
%}~~~~\text{quantum}
\]
Drinfel'd \cite{Dri} and Jimbo \cite{Jim} have shown how to
define a quantization $\uqg$ of the enveloping algebra $U(\lie)$ of
any simple complex Lie algebra $\lie$. These quantized enveloping
algebras have proven to be important in several branches of 
mathematics and physics and have been studied in detail. In contrast,
very little is know about the quantization of the Lie algebras
themselves.

The approach to the quantization of Lie algebras which was initiated
in \cite{qliea,qlieb} consists of making use of the known quantization
of the corresponding enveloping algebras. The quantum Lie algebras
are extracted from the quantized enveloping algebras in the same
way in which the unquantized Lie algebras can be extracted from the
unquantized enveloping algebras.

In an alternative approach, initiated by Woronowicz \cite{Wor}, one
can extract quantum Lie algebras from the formalism of bicovariant
differential calculi on quantum groups. This approach always leads
to quantum Lie algebras which have a larger dimension than their
classical counterpart. We now know that these algebras are not simple
and that our quantum Lie algebras are the simple subalgebras
of these with the correct dimension. Here we will deal only with 
the simple quantum Lie algebras
obtained from our the algebraic approach. 
For works on the Woronowicz algebras see e.g.
\cite{others}. For an approach specific to quantum $\sln$ but
related to our general approach
see \cite{sud}.

In this paper we will use the simplest example, namely the three-dimensional
Lie algebra $\slt$ so familiar to physicists, to introduce the general
results about quantum Lie algebras which have recently been obtained.
We will also give a matrix representation of
this algebra which has not yet been published.

\section{The Lie algebra $\slt$}\label{sectsl2}

Complex Lie algebras $\lie$ in general are vector spaces over $\CC$ equipped
with a non-associative product, commonly denoted as the Lie bracket. This
is a linear map $\lb:\lie\ot\lie\rightarrow\lie$ which satisfies
\begin{align}
[a,b]&=-[b,a]&&\text{antisymmetry,}\label{antisym}\\
[a,[b,c]]&=[[a,b],c]+[b,[a,c]]&&\text{Jacobi identity.}\label{jacobi}
\end{align}
The complex simple Lie algebra $\slt$ is spanned as a vector space
by three elements $X^+,X^-$ and $H$. The Lie bracket
is given by
\begin{equation}\label{lierels}
[X^+,X^-]=H,~~~~~~
[H,X^\pm]=\pm 2 X^\pm.
\end{equation}
Together with the antisymmetry property and the bilinearity these three
relations define
the Lie bracket on the whole algebra uniquely.

\section{The enveloping algebra $U(\slt)$\label{sectU}}

The enveloping algebra $U(\slt)$ is the associative unital algebra over
$\CC$
generated by the three generators $X^+,X^-$ and $H$ and the commutation
relations
\begin{equation}\label{commrels}
X^+ X^- - X^- X^+ =H,~~~~
H X^\pm - X^\pm H =\pm 2 X^\pm.
\end{equation}
In other words: $U(\slt)$ contains all possible ordered polynomials in the
three generators but two such polynomials are equal if they are related
by the above commutation relations. It can be seen that
the relations allow one to commute all $X^-$ to the left and all 
$X^+$ to the right. Thus as a basis for $U(\slt)$ one can choose
$\{(X^-)^n H^m (X^+)^l|n,m,l\in\NN\}$, known as the 
Poincar\'e-Birkhoff-Witt basis.

The enveloping algebra $U(\slt)$ is clearly infinite dimensional. 
It contains $\slt$ as the subspace spanned by $X^-,X^+$ and $H$.
This subspace is closed under the commutator and the commutator
coincides with the Lie bracket as defined in \eqref{lierels}. Because
of this the mind of a physicist tends not to distinguish between
Lie brackets and commutators. Below however it will be crucial to
keep the two concepts clearly separated.

The enveloping algebra is a Hopf algebra. In this paper we will only
need to know that this implies that one can define an action of the
enveloping algebra on itself, the so called adjoint action. For
$U(\slt)$ it is defined by
\begin{equation}
\ad{X^\pm}a=X^\pm a - a X^\pm,~~~~~~
\ad{H}a=H a - a H,~~~~~\forall a\in U(\slt).
\end{equation}
Thus for the generators the adjoint action is just the commutator.
For products of generators the adjoint action is obtained
from the above by the defining property of an action, i.e.,
$\ad{ab}=\ad{a}\ad{b}$.
\section{The quantized enveloping algebra $U_h(\slt)$}

The quantized enveloping algebra $U_h(\slt)$ \cite{Dri,Jim,Andrew} too is an 
associative unital algebra
generated by the three generators $X^+,X^-$ and $H$. 
However it is an algebra over $\ch$, the ring of formal power series in
an indeterminate $h$ (which in physical applications may be related
to Planck's constant, but not necessarily linearly so). The
commutation relations are deformed with respect to \eqref{commrels}.
They now read
\begin{equation}\label{qcommrels}
X^+ X^- - X^- X^+ =\frac{q^H-q^{-H}}{q-q^{-1}},~~~~
H X^\pm - X^\pm H =\pm 2 X^\pm,
\end{equation}
where $q=e^h$. Thus the commutator of $X^+$ and $X^-$ gives an
infinite power series in h. The first term in the series is just $H$,
as classically, but the higher order terms in $h$ (the quantum corrections)
are non-linear in $H$. 

The important property of the deformation \eqref{qcommrels} is that it
still defines a Hopf algebra. However also the Hopf algebra structure
is deformed and this leads in particular to a deformed adjoint action
\begin{align}\label{qadj}
\ad{X^\pm}a&=X^\pm a q^{H/2}- q^{\mp 1} q^{H/2} a X^\pm,\nonumber\\
\ad{H}a&=H a - a H,&&\forall a\in U(\slt).
\end{align}

\section{A quantum Lie algebra $\lqslt$ inside $U_h(\slt)$\label{sectlqslt}}

As explained in \secref{sectU} the Lie algebra $\slt$
can be viewed as a subspace of the enveloping algebra $U(\slt)$
which is spanned by the generators $X^+,X^-$ and $H$ and on
which the Lie bracket is given by the commutator. We would now
like to obtain the quantum Lie algebra $\lqslt$ in a similar manner
from the quantized enveloping algebra $U_q(\slt)$. 

However, in the quantum case the space\footnote{To be pedantic, 
because $\ch$ is a ring
and not a field, we should not speak of vector spaces but rather of
$\ch$-modules. However in this paper we would prefer not to dwell
on such subtleties.} spanned by the generators $X^+,X^-$ and $H$
does not close under the commutator. The non-linear terms in H in
the commutation relations \eqref{qcommrels} create a problem.
The first idea is to
replace the role of the commutator by the adjoint action. 
As we had
seen the two coincide in the classical case but differ in the quantum
case.
If we find a three-dimensional subspace of $U_q(\slt)$ which is
closed under the adjoint action then we can define a quantum
Lie bracket on this space by
\begin{equation}\label{qbdef}
[a,b]_h\equiv\ad{a}b.
\end{equation}
The space spanned by
\begin{equation}\label{qliegen}
\xpm=\sqrt{\tfrac{2}{q+q^{-1}}}\,q^{-H/2} X^\pm,~~~~~~
\hh=\tfrac{2}{q+q^{-1}}\left(q X^+ X^--q^{-1} X^- X^+\right),
\end{equation}
satisfies this requirement. Indeed, using \eqref{qadj},
one can calculate the adjoint
action of these elements on each other. The reader is urged to
perform these calculations. He will find
\begin{align}\label{qlierels}
[\xp,\xm]_h&=\hh,&
[\xm,\xp]_h&=-\hh,\nonumber\\
[\hh,\xpm]_h&=\pm 2 q^{\pm 1} \xpm,&
[\xpm,\hh]_h&=\mp 2 q^{\mp 1} \xpm\nonumber\\
[\hh,\hh]_h&=2(q-q^{-1}) \hh,&
[\xpm,\xpm]_h&=0.
\end{align}
These quantum Lie bracket relations are the quantum analog of the
Lie bracket relations \eqref{lierels}. To zeroth order in $h$
they agree with the classical Lie bracket relations.

\section{$q$-antisymmetry}

We stress that the quantum Lie algebra $\lqslt$ is not a Lie algebra.
The quantum Lie bracket defined by \eqref{qlierels} does not satisfy
the properties of antisymmetry \eqref{antisym} and Jacobi \eqref{jacobi}.
Instead it satsifies an interesting generalization of antisymmetry which
involves the operation $q\rightarrow 1/q$. The details are as follows:

We define $q$-conjugation $\qconj:\ch\rightarrow\ch$ as the
$\CC$-linear ring automorphism defined by $\t{h}=-h$ (and thus
$\t{q}=1/q$). We extend this to a $q$-conjugation on $\lqslt$ by
defining 
\begin{equation}
\left(a\, \xp+b\,\xm+c\,\hh\right)^\qconj=\t{a}\,\xp+\t{b}\,\xm+\t{c}\,\hh.
\end{equation}
Then the quantum Lie bracket satisfies
\begin{equation}\label{qantisym}
[a,b]_h^\qconj=-[\t{b},\t{a}]_h.
\end{equation}
We call this property $q$-antisymmetry. 
This property can easily be verified for the quantum Lie bracket relations
\eqref{qlierels} even though it is not at all evident from the definition
\eqref{qbdef} of the quantum Lie bracket. As shown in \cite{qliec} all
quantum Lie algebras (defined below) possess this $q$-antisymmetry.

We have not yet discovered the $q$-analog of the Jacobi identity.

\section{General definition of $\lqlie$}

The quantum Lie algebra $\lqslt$ constructed in \secref{sectlqslt}
is the simplest example of the following general definition.

\begin{definition}\label{defqlie}
A \dem{quantum Lie algebra $\lqlie$
inside $\uqg$} is a finite-dimensional 
irreducible $\text{ad}$ - submodule
of $\uqg$ endowed with the 
\dem{quantum Lie bracket} 
$[a,b]_h{}=\ad{a}{b}$
such that
\begin{enumerate}
\vspace{-1mm}
\item $\lqlie$ is a deformation of $\lie$, i.e., there is an
algebra isomorphism
$\lqlie\cong\lie\pmod h$.
\vspace{-1mm}
\item
$\lqlie$ is invariant under the $q$-Cartan involution $\t{\theta}$, 
the $q$-antipode $\t{S}$ and any diagram automorphism $\tau$ of $U_q(\lie)$.
\end{enumerate}
\vspace{-1mm}
%A \dem{weak quantum Lie algebra $\wqlie$} is defined similarly but
%without the requirement 2.
\end{definition}
Property 2 plays an important role in the investigations into
the general structure of quantum Lie algebras in \cite{qliea}. In particular
it allows the definition of a quantum Killing form. We refer the
reader to the paper \cite{qliea} for more information on these matters.
It is shown in \cite{qliea} that given any module satisfying all
properties of the definition except property 2 one can always
construct from it a quantum Lie
algebra $\lqlie$ which satisfies property 2 as well. Thus this
extra requirement is not too strong. 

While finding a quantum Lie algebra $\lqslt$ inside $U_h(\slt)$ was
easy, performing the similar task for other algebras is much more
involved. However, as reported in \cite{qliea}, it has been done
by using the computer algebra program Mathematica for the Lie
algebras $\lie=\mathfrak{sl}_3,\mathfrak{sl}_4,\mathfrak{so}_5=
\mathfrak{sp}_4$ and $G_2$. There is also a method for constructing
quantum Lie algebras in general using the universal R-matrix. This
method has been applied in \cite{qlieb} to obtain quantum Lie algebras
for $\lie=\sln$ for all $n$. The method is described also in
\cite{qliec}. This construction proves in particular that quantum Lie
algebras exist for all $\lie$. However the determination of the structure
of these quantum Lie algebras has not yet been performed except in the
above mentioned cases.

\section{The abstract quantum Lie algebra $\qslt$}

The space spanned by the three generators $\xp,\xm$ and $\hh$ given
in \eqref{qliegen} is not the only three dimensional subspace of 
$U_h(\slt)$ which is closed under the adjoint action and gives rise
to a quantum Lie algebra $\lqslt$ according to the general definition
in the previous section. 
But it is easy to convince oneself
of the fact that any quantum
Lie algebra $\lqslt$ inside $U_h(\slt)$ is spanned by three
elements of the form
\begin{equation}\label{qliegenc}
\xpm=\sqrt{\tfrac{2}{q+q^{-1}}}\,q^{-H/2} X^\pm \poly{C},~~~~~~
\hh=\tfrac{2}{q+q^{-1}}\left(q X^+ X^--q^{-1} X^- X^+\right)\poly{C},
\end{equation}
where $\poly{C}$ can be any polynomial in the Casimir element
\begin{equation}\label{casimir}
C=\frac{1}{(q^3+q^{-3})}\left(
(q-q^{-1})^2 X^+ X^-+q^{H-1}+q^{-H+1}\right)
\end{equation}
whose coefficients sum to $1$.
Because of the properties of the Casimir element, all these
quantum Lie algebras lead to the same quantum Lie bracket
relations \eqref{qlierels}. Thus all quantum
Lie algebras $\lqslt$ are isomorphic as algebras. This
defines an abstract quantum Lie algebra $\qslt$. $\qslt$ is
the algebra spanned by three abstract generators $\xp,\xm$ and $\hh$
with the Lie bracket relations \eqref{qlierels}.  The concrete
quantum Lie algebras $\lqslt$ are just
different embeddings of $\qslt$ into $U_h(\slt)$. 

Similarly, as proven in \cite{qliec}, 
there is a unique abstract quantum Lie algebra $\qlie$
for any simple complex Lie algebra $\lie$.
\begin{theorem}
All concrete quantum Lie algebras $\lqlie$ for the same $\lie$
are isomorphic to a unique abstract quantum Lie algebra $\qlie$.
\end{theorem}
Furthermore it has been shown that the structure constants of
$\qlie$ are equal to the inverse $q$-Clebsch-Gordon coefficients for
adjoint$\times$adjoint into adjoint.

\section{Representations of $\qslt$\label{sectreps}}

An $n$-dimensional representation of a Lie algebra is a linear map
$\pi$ from the Lie algebra into the $n\times n$ matrices,
$\pi:\lie\rightarrow\mat_n(\CC)$, such that the Lie bracket is
realized as the matrix commutator, i.e.,
\begin{equation}
\pi([a,b])=\pi(a)\pi(b)-\pi(b)\pi(a).
\end{equation}
Such maps $\pi$ exist because the commutator also possesses the 
defining properties of the Lie bracket, namely antisymmetry \eqref{antisym}
and Jacobi \eqref{jacobi}.

To arrive at a good definition of a representation of a quantum Lie
algebra we have to find a $q$-generalization of the commutator which
is $q$-antisymmetric in the sense of \eqref{qantisym}. We propose
the following definition:

\begin{definition}
An $n$-dimensional \dem{representation of a quantum Lie algebra} $\qlie$
consists of
\vspace{-2mm}
\begin{itemize}
\item a linear map $\pi:\qlie\rightarrow\mat_n(\ch)$, 
%$\pi:\qlie\rightarrow\mat_n(\ch)$, and of
\vspace{-2mm}
\item a $q$-conjugation $\qconj:\mat_n(\ch)\rightarrow\mat_n(\ch)$,
\end{itemize}
\vspace{-1mm}
such that
\begin{align}
\vspace{-2mm}
\pi(\t{a})&=\t{\pi}(a),\\
\pi([a,b]_h)&=\pi(a)\pi(b)-\left(\t{\pi}(b)\t{\pi}(a)\right)^\qconj.
\label{qcomm}
\end{align}
\end{definition}
By definition \cite{qliea} a $q$-conjugation is $q$-linear, i.e., 
$(\lambda\,a)^\qconj=\t{\lambda}\,\t{a}~
\forall\lambda\in\ch,a\in\mat_n(\ch)$, and is an involution, i.e.,
$\t{\t{a}}=a~\forall a\in\mat_n(\ch)$.

To illustrate this definition we will give the 2-dimensional
representation of $\qslt$. The representation matrices are
\begin{align}
\pi(\xp)&=\sqrt{\tfrac{q+q^{-1}}{2}}
\begin{pmatrix}0&1\\0&0\end{pmatrix},&
\pi(\hh)&=\begin{pmatrix}q&0\\0&-q^{-1}\end{pmatrix},&
\pi(\xm)=\pi(\xp)^t
\end{align}
and the $q$-conjugation is given by
\begin{align}\label{repmat}
\begin{pmatrix}1&0\\0&0\end{pmatrix}^\qconj &=
\frac{1}{q+q^{-1}}\begin{pmatrix}2q&0\\0&q-q^{-1}\end{pmatrix},&
\begin{pmatrix}0&1\\0&0\end{pmatrix}^\qconj&=
\begin{pmatrix}0&1\\0&0\end{pmatrix},\nonumber\\[4mm]
\begin{pmatrix}0&0\\0&1\end{pmatrix}^\qconj &=
\frac{1}{q+q^{-1}}\begin{pmatrix}q^{-1}-q&0\\0&2q^{-1}\end{pmatrix},&
\begin{pmatrix}0&0\\1&0\end{pmatrix}^\qconj&=
\begin{pmatrix}0&0\\1&0\end{pmatrix}.
\end{align}
The reader is urged to check that the $q$-commutators of the representation
matrices in \eqref{repmat} do indeed reproduce the algebra \eqref{qlierels}.

The author has a construction for representations of $\qslt$ of any
dimension.

\section{Discussion}

We have reviewed quantum Lie algebras by using the explicit example
of $\slt$. Quantum Lie algebras have been defined as certain
subspaces $\lqlie$ of the quantized enveloping algebras $\uqg$. It has
been found that as algebras all the $\lqlie$ are isomorphic to an 
abstract quantum Lie algebra $\qlie$. We have seen that $q$-conjugation 
$q\rightarrow 1/q$ plays
a central role in the theory of quantum Lie algebras. In particular,
the quantum Lie bracket has turned out to be $q$-antisymmetric in the sense
of \eqref{qantisym}. This has lead to the definition \eqref{qcomm} of
a $q$-commutator to represent the quantum Lie bracket.

Many definitions for $q$-commutators can be found in the literature.
They generally differ from the usual commutator by multipying the
terms by certain powers of $q$. Our definition \eqref{qcomm}
of the $q$-commutator is quite different and arises naturally in
the theory of quantum Lie algebras. It will be interesting to
study its physical applications.

Drinfel'd has introduced a quantized enveloping algebra $\uqg$ for any
complex simple Lie algebra $\lie$. There is no definition of what a
quantized enveloping algebra is in general. Similarly we have so far
defined the concept of a quantum Lie algebra $\qlie$ 
only for complex simple (finite-dimensional) Lie algebras $\lie$.
There are however indications that an axiomatic definition of quantum
Lie algebras can be obtained which parallels the axiomatic definition of 
Lie algebras
through the properties of antisymmetry \eqref{antisym} and Jacobi
\eqref{jacobi}.
Antisymmetry and Jacobi identity are the necessary and sufficient conditions
for an algebra to have a representation in terms of commutators. 
Similarly $q$-antisymmetry and
$q$-Jacobi should be the necessary and sufficient conditions for
an algebra to have a representation in terms of $q$-commutators. 
The $q$-antisymmetry \eqref{qantisym} is clearly a necessary
condition but we are still searching for the $q$-Jacobi identity which
gives a sufficient condition.
It could be hoped that this will
then also finally lead to an axiomatic definition of quantized 
enveloping algebras.

For more information on quantum Lie algebras visit their WWW site at

http://www.mth.kcl.ac.uk/\~{}delius/q-lie.html.

\end{document}